\documentclass[graybox]{svmult}

\usepackage{type1cm}        
\usepackage{makeidx}         
\usepackage{graphicx}        
\usepackage{multicol}        
\usepackage[bottom]{footmisc}

\usepackage{newtxtext}       %
\usepackage{newtxmath}       
\usepackage{amsmath}

\newcommand{\F}{\mathcal{F}}
\newcommand{\Sc}{\mathcal{S}}
\newcommand{\A}{\mathcal{A}}
\newcommand{\R}{\mathbb{R}}
\newcommand{\N}{\mathbb{N}}

\allowdisplaybreaks
\makeindex             

\begin{document}

\title*{Linear Equations in the Ring of $\mathcal{S}\left(\mathcal{A}\right)$-Linearly Correlated Fuzzy Numbers}
\titlerunning{Linear Equations in the Ring of $\mathcal{S}\left(\mathcal{A}\right)$-Linearly Correlated Fuzzy Numbers} 
\author{Beatriz Laiate and Peter Sussner}
\institute{Beatriz Laiate \at School of Applied Mathematics - Getúlio Vargas Foundation (FGV/EMAp), Rio de Janeiro, Brazil \email{beatrizlaiate@gmail.com, beatriz.laiate@fgv.br}
\and Peter Sussner \at Institute of Mathematics, Statistics and Scientific Computing, Universidade Estadual de Campinas (IMECC/Unicamp), Campinas, Brazil \email{sussner@unicamp.br}}
%
%
\maketitle

\abstract*{This paper investigates the solutions of a family of certain linear fuzzy arithmetic equations
that involve fuzzy numbers belonging to certain finite-dimensional vector spaces of $\mathbb{R}_{\mathcal{F}}$, called  $\mathcal{S}\left(\mathcal{A}\right)$-linearly correlated fuzzy numbers. Here, $\A$ stands for a strongly linearly independent (SLI) set of fuzzy numbers.
The arithmetic operations in the aforementioned linear equations
are the sum in the vector space $\Sc(\A)$ and the so-called $\psi$-cross product
that turn  $\Sc(\A)$ into a commutative ring.}


\abstract{This paper investigates the solutions of a family of certain linear fuzzy arithmetic equations
that involve fuzzy numbers belonging to certain finite-dimensional vector spaces of $\mathbb{R}_{\mathcal{F}}$, called  $\Sc(\A)$-linearly correlated fuzzy numbers. Here, $\A$ stands for a strongly linearly independent (SLI) set of fuzzy numbers.
The arithmetic operations in the aforementioned linear equations
are the sum in the vector space $\Sc(\A)$ and the so-called $\psi-$cross product
that turn $\Sc(\A)$ into a commutative ring.}
\keywords{Fuzzy numbers, Fuzzy arithmetic operations, Fuzzy linear equation, Ring of $\mathcal{S}\left(\mathcal{A}\right)$-linearly correlated fuzzy Numbers.}

\section{Introduction} \label{sec:1}

Fuzzy arithmetic equations can be seen as the first step to solving fully fuzzy linear systems, for which matrices with fuzzy entries are considered. As well as a proper interpretation of its solution \cite{lodwick2015interval}, fuzzy arithmetic operations are essential to this study.

The choice of the definitions of arithmetic operations for fuzzy numbers determines how to deal with fuzzy equations. Zadeh's extension principle yields the first approach.  It is a well-known fact that the usual arithmetic operations for fuzzy numbers  $A, B \in \mathbb{R}_{\mathcal{F}}$ are determined by Zadeh's extension of the arithmetic operations on real numbers, that is, $A \otimes B = \hat{\otimes}\left(A,B\right)$, where $\otimes \in \left\{+,-,\times,\div\right\}$. In particular, the usual product operation $A$ and $B$, $A \times B$, is defined via Zadeh's Extension of the function $f(x,y)=xy$ applied to $(A,B)$ \cite{zadeh1975concept}.

An alternative approach toward fuzzy arithmetic operations was presented by Dubois and Prade \cite{dubois1978operations,dubois1983inverse,dubois1993fuzzy} who introduced the notion of LR fuzzy numbers. Fuzzy linear equations have been studied by several researchers  \cite{behera2014solving,biacino1989equations,buckley1990solving,esmi2019some,ghanbari2022new,sanchez1984solution,sevastjanov2009new}. 
Usually, a fuzzy number does not have an inverse with respect to an arithmetic operation 
 \cite{yager1980lack}. 

This paper deals with fuzzy arithmetic operations defined in finite-dimensional vector spaces of fuzzy numbers. 
The focus of this paper is on solving linear equations of the form
\begin{equation} \label{eq:main}
    A \odot_\psi X +_\psi B = C,    
\end{equation}
where $A,B,C$ are linear combinations of elements of an SLI set of fuzzy numbers $\left\{1 = A_1, A_2, \ldots,A_n\right\} \subset \mathbb{R}_{\mathcal{F}}^\wedge \subset \R_F $ \cite{esmi2021banach,laiate2021cross}. Here, $\R_\F$ denotes the class of fuzzy numbers and $\mathbb{R}_{\mathcal{F}}^\wedge$ stands for the set of all $A \in \R_\F$ whose core has only one element. Moreover, $\odot_\psi$ and $+_\psi$ represent the $\psi-$cross product in  $\Sc(\A)$ and the sum that is induced by the isomorphism between  $\R^{n}$ and $\Sc(\A)$,
respectively. 
As mentioned before, $(\Sc(\A), +_\psi, \odot_\psi)$ is a commutative ring and Eq. \eqref{eq:main} must be
solved in this algebraic structure. Since $A \in \Sc(\A) \subset \R_{\F}^\wedge$ has an inverse with respect to $\odot_\psi$ if and only if $[A]_1 \neq \{0\}$, the solution of the general linear equation above can be 
determined under some weak conditions. 

The remainder of this paper is organized as follows: Section 2 reviews some preliminary concepts. Section 3 extends the results presented in \cite{laiate2021bidimensional} and introduces a linear fuzzy-valued function in the space $\mathcal{S}\left(\mathcal{A}\right)$, proving that it is a bijection. Section $4$ uses the previous results to determine the solution of Eq. \ref{eq:main} under some weak conditions. We finish with some concluding remarks in Section 5.

\section{Preliminaries} \label{sec:2}

A fuzzy subset $A$ of a universe $X$, aka a fuzzy set on $X$, can be identified with its membership function $\mu_A: X \to  [0,1]$. Thus, it makes sense to write $A(x)$ instead of 
 $\mu_A(x)$ for any $x \in X$.
The class of fuzzy sets on $X$ is denoted using the symbol
${\mathcal{F}}(X)$. Recall that every $A \in {\mathcal{F}}(X)$ is determined by its $\alpha$-levels, aka $\alpha$-cuts, $[A]_\alpha$. For every $\alpha \in (0,1]$, we have that 
\begin{equation}
[A]_\alpha = \left\{x \in X \, | \, A(x) \geq \alpha \right\}. 
\end{equation}
The {\it support} of $A$, denoted $\mbox{supp}(A)$, is given by  $\mbox{supp}(A) = \{x \in X \, | \, A(x) > 0 \}$.
A fuzzy set $A$ on $\mathbb{R}$  is called a {\it fuzzy number} if $A$ satisfies the following conditions:
\begin{enumerate}
    \item $A$ is normal, i.e., $[A]^1 \neq \emptyset$;
    \item  $A$ is fuzzy convex;
    \item  $\mbox{supp}(A)$  is bounded;
    \item $\mu_A$ is upper semicontinuous.
\end{enumerate}
The symbol $\mathbb{R}_{\mathcal{F}}$ denotes the class of fuzzy numbers. For a fuzzy number $A$, the $0-$level of $A$ can be defined 
as the closure of its support. With this definition, a fuzzy number $A$ can simply be characterized as an element of ${\mathcal{F}}(\mathbb{R})$ whose
$\alpha$-cuts, where $\alpha \in [0,1]$,  are non-empty, compact intervals.
Thus, if $A \in \mathbb{R}_{\mathcal{F}}$, then  $\left[A\right]_{\alpha}$ is of the form $\left[a_{\alpha}^-,a_{\alpha}^+\right]$, where $a_{\alpha}^-$ and $a_{\alpha}^+ \in \mathbb{R}$ are the lower and upper endpoints of $\left[A\right]_{\alpha}$, respectively. 
For each $\alpha \in [0,1]$, the length of $\left[A\right]_\alpha$ is defined as len$\left(\left[A\right]_\alpha\right) = a_{\alpha}^+ - a_{\alpha}^-$. In particular, the length of $\left[A\right]_0$ is called the {\it diameter} of $A$ and denoted diam$(A)$. The diameter of a fuzzy number can be seen as a measure of the uncertainty, or {\it fuzziness}, of a quantity represented by $A \in \mathbb{R}_{\mathcal{F}}$.

A fuzzy number $A \in \mathbb{R}_{\mathcal{F}}$ is said to be {\it symmetric with respect to $x \in \mathbb{R}$}, expressed symbolically as $\left(A|x\right)$,  if $A(x-y)=A(x+y)$ for all $y \in \mathbb{R}$. If there exist  $x \in \mathbb{R}$ such that $A$ is symmetric  with respect to $x$, then
$A$ is said to be {\em symmetric}. Otherwise, $A$ is called  non-symmetric. The Hausdorff-Pompeiu distance between two fuzzy numbers is given by
\begin{equation} \label{eq:heusdorffdist}
    \mathcal{D}_\infty \left(A,B\right) = \sup_{\alpha \in [0,1]} \max\left\{|a_{\alpha}^--b_{\alpha}^-|,|a_{\alpha}^+-b_{\alpha}^+|\right\},
\end{equation}
\noindent for all $A,B \in \mathbb{R}_{\mathcal{F}}$ given levelwise by $\left[A\right]_{\alpha} = \left[a_{\alpha}^-,a_{\alpha}^+\right]$ and $\left[B\right]_{\alpha} = \left[b_{\alpha}^-,b_{\alpha}^+\right]$, $\forall \alpha \in [0,1]$. Recall that $\left(\mathbb{R}_{\mathcal{F}},\mathcal{D}_\infty\right)$ is a complete, but not a separable metric space \cite{diamond1990metric}. 

\begin{definition} \label{def:zadehextension} \cite{nguyen1978note,zadeh1975concept} Let $X, Y \neq \emptyset$,  $A \in \mathcal{F}(X)$, and $f:X \to Y$. The Zadeh extension of $f$ at $A$, denoted by $\hat{f}(A)$, is given by 
    \begin{equation} \label{eq:zadehextension}
        \hat{f}\left(A\right)(y) = \begin{cases}
            \begin{array}{cl} \displaystyle
                \sup_{x \in f^{-1}(y)} A(x) & \; \mbox{if} \; f^{-1}\left(y\right) \neq \emptyset\\
                0 & \; \mbox{if} \; f^{-1}\left(y\right) = \emptyset,
            \end{array}
        \end{cases},
    \end{equation}
    \noindent where $f^{-1}\left(y\right) = \left\{x \in X : f(x)=y\right\}$, $\forall y \in Y$.
\end{definition}

Using the notion of the Cartesian product of fuzzy sets,  a function of multiple variables of the form $f: X_1 \times X_2 \ldots X_n \to Y$ can be extended
to a mapping  $\hat{f}: \F(X_1) \times \F(X_2) \ldots \F(X_n) \to \F(Y)$. Recall that 
$A_1 \times A_2 \times \ldots \times A_n$, where $A_i \in \F(X_i)$ for $i = 1, 2, \ldots n$, is given by
\begin{equation}
(A_1 \times A_2 \times \ldots \times A_n)(x_1, x_2, \ldots, x_n) = \min\left\{A_1(x_1),\ldots,A_n(x_n)\right\}.
 \end{equation}
The mapping  $\hat{f}: \F(X_1 \times X_2 \ldots X_n) \to \F(Y)$ induces a mapping
$\F(X_1) \times \F(X_2) \ldots \F(X_n) \to \F(Y)$ that is also denoted  $\hat{f}$ via 
the equation
\begin{equation}
\hat{f}(A_1, A_2, \ldots, A_n) = \hat{f}(A_1 \times A_2 \ldots A_n) \, \, \forall A_i \in \F(X_i), \,  i = 1, 2, \ldots, n.
 \end{equation}
Hence, Zadeh's extension of $f$ at $(A_1,\ldots,A_n)$ is given by

\begin{equation} \label{eq:zadehextension2}
    \hat{f}(A_1,\ldots,A_n)(y) = \begin{cases}
    \begin{array}{cl}
    \displaystyle \sup_{(x_1,\ldots,x_n) \in f^{-1}(y)} \min\left\{A_1(x_1),\ldots,A_n(x_n)\right\} \quad &\mbox{if} \; f^{-1}(y) \neq \emptyset\\
    0     &\mbox{if} \; f^{-1}(y) = \emptyset
    \end{array}       
    \end{cases},
\end{equation}
\noindent where $f^{-1}\left(y\right) = \left\{(x_1,\ldots,x_n) \in X_1 \times \ldots \times X_n : f(x_1,\ldots,x_n)=y\right\}$. \cite{bede,fuller1998fuzzy}

Under the hypothesis of continuity of $f$, Zadeh's Extension of $f$ applied to some fuzzy number $A \in \mathbb{R}_{\mathcal{F}}$ can be obtained levelwise, as stated in the following \cite{wang09}:

\begin{theorem} \label{the:zadehsextension} 
    If a real function $f$ of $n$ real variables is continuous and
if $A_1, A_2, \ldots, A_n$ are fuzzy numbers, then
    \begin{equation}
        \left[\hat{f}\left(A\right)\right]_{\alpha} = f\left(\left[A\right]_{\alpha}\right), \quad \forall \alpha \in (0,1].
    \end{equation}
\end{theorem}

Theorem \ref{the:zadehsextension} assures that Zadeh's Extension applied to a fuzzy number $A \in \mathbb{R}_{\mathcal{F}}$ of a continuous vector-valued function is completely determined by the function evaluated in the $\alpha-$cuts of $A$. This paper focuses on a special class of fuzzy numbers, namely, those contained in finite-dimensional vector spaces of $\mathbb{R}_{\mathcal{F}}$. In the following, we recall the notions introduced in previous papers \cite{esmi2021banach,esmi2022calculus}.

Let $\mathcal{A} = \left\{A_1,\ldots,A_n\right\} \subset \mathbb{R}_{\mathcal{F}}$ be a finite set of fuzzy numbers. The set of all linear combinations of $A_1, \ldots, A_n$ is denoted by
\begin{equation} \label{eq:linearcomb}
    \mathcal{S}\left(\mathcal{A}\right) = \left\{q_1A_1+\ldots+q_nA_n : q_i \in \mathbb{R}, i=1,\ldots, n\right\},
\end{equation}
\noindent where $'+'$ and $q_iA_1$ stand for the usual operations of sum and scalar multiplication in $\mathbb{R}_{\mathcal{F}}$. If $B \in \mathcal{S}\left(\mathcal{A}\right)$, then there exist $q_1,\ldots,q_n \in \mathbb{R}$ such that $B = q_1A_1+\ldots+q_nA_n$ or, levelwise,
\begin{equation*}
    \left[B\right]_{\alpha} = q_1\left[A_1\right]_{\alpha} + \ldots q_n\left[A_n\right]_{\alpha}, \quad \forall \alpha \in [0,1].
\end{equation*}

The next definition recalls the notion of  {\it Strong Linear Independence} of a set of fuzzy numbers.

\begin{definition} \label{def:SLIsets}
    A set of fuzzy number $\mathcal{A} = \left\{A_1,\ldots,A_n\right\}$ is said to be strongly linearly independent (SLI, for short) if we have for every $B = q_1A_1 + \ldots + q_nA_n \in \mathcal{S}\left(\mathcal{A}\right)$ that 
    \begin{equation} \label{eq:SLI}
        \left(B | 0\right) \Rightarrow q_1 = \ldots = q_n = 0.
    \end{equation}
\end{definition}

Definition \eqref{def:SLIsets} resembles the notion of {\it linear independence} of a set of vectors in a vector space. Recall $\mathbb{R}_{\mathcal{F}}$ together with the usual arithmetic operations is a quasilinear space and that strong linear independence can also be defined for infinite sets of fuzzy numbers \cite{esmi2021banach}.
The next theorem provides a complete characterization of SLI sets of fuzzy numbers.

\begin{theorem} \label{the:SLIisomorphism}
    The set $\mathcal{A} = \left\{A_1,\ldots,A_n\right\} \subset \mathbb{R}_{\mathcal{F}}$ is SLI if, and only if, {\color{blue}$\Sc(\A)$} is a real vector space that is isomorphic to $\mathbb{R}^n$ under the isomorphism $\psi: \mathbb{R}^n \rightarrow \mathcal{S}\left(\mathcal{A}\right)$ given by
    \begin{equation} \label{eq:SLIisomorphism}
        \psi\left(x_1,\ldots,x_n\right) = x_1A_1 + \ldots + x_nA_n.
    \end{equation}
\end{theorem}

If $B \in \mathcal{S}\left(\mathcal{A}\right) $ for some SLI set $\mathcal{A} \subset \mathbb{R}_{\mathcal{F}} \left\{A_1,\ldots,A_n\right\}$, then $B$
is said to be an {\it $\mathcal{S}\left(\mathcal{A}\right)-$linearly correlated fuzzy number}, also called
$\mathcal{S}-$linearly correlated fuzzy number. (Since this concept depends on an SLI set $\A$, we 
prefer to speak of the former instead of the latter.)
In this case, 
there exists a unique vector  $(q_1,\ldots,q_n) \in \mathbb{R}^n$ such that $B = \psi\left(q_1,\ldots,q_n\right) = q_1A_1+\ldots+q_nA_n$. 

From now on, let us  assume that $\mathcal{A}=\left\{1,A_2,\ldots,A_n\right\}$, where $A_i \in \mathbb{R}_{\mathcal{F}} \backslash \mathbb{R}$, $i =2, \ldots,n$ and that  
$$\mathcal{A} \subset \mathbb{R}_{\mathcal{F}}^{\wedge}=\left\{A \in \mathbb{R}_{\mathcal{F}} : \left[A\right]_1 \mbox{has only one element}\right\} \subset \mathbb{R}_{\mathcal{F}}.$$ 
For this choice of $\A$, an  $\mathcal{S}(\A)-$linearly correlated fuzzy number $B = q_1 + q_2A_2 + \ldots + q_nA_n \in \mathcal{S}\left(\mathcal{A}\right)$ represents the sum of the crisp quantity $q_1$ and the  fuzzy quantities $q_2A_2, \ldots, q_nA_n$.

The arithmetic operations in $\mathcal{S}\left(\mathcal{A}\right)$ arise from the bijection given by Eq. \eqref{eq:SLIisomorphism}, as stated in the following. For all $\lambda \in \mathbb{R}$ and $B, C \in \mathcal{S}\left(\mathcal{A}\right)$ given by $B = q_1 + q_2A_2 + \ldots + q_nA_n$ and $C = p_1 + p_2A_2 + \ldots p_nA_n$, respectively, we have
\begin{equation} \label{eq:sumpsi}
    B +_{\psi} C \, \dot{=} \, (q_1+p_1) + (q_2+p_2)A_2 + \ldots + (q_n+p_n)A_n    
\end{equation}
and 
\begin{equation} \label{eq:prodpsi}
    \lambda \cdot_\psi B \, \dot{=} \, (\lambda q_1) + (\lambda q_2)A_2 + \ldots + (\lambda q_n)A_n.
\end{equation} 

Eqs. \eqref{eq:sumpsi} and \eqref{eq:prodpsi} imply that $$B -_{\psi} C \, \dot{=} \, (q_1-p_1) + (q_2-p_2)A_2 + \ldots + (q_n-p_n)A_n.$$ 
The product of $B$ and $C$ is given by
$$B \odot_\psi C = cB +_\psi bC -_\psi bc,$$ where $\left[B\right]_1 = \left\{b\right\}$ and $\left[C\right]_1 = \left\{c\right\}$. If $c \neq 0$, then $B \div_{\psi} C \, \dot{=} \, B \odot_\psi C_{\psi}^{-1}$, where 

$$C_{\psi}^{-1} \, \dot{=} \, \left(\frac{2}{c}-\frac{p_1}{c^2}\right) - \frac{p_2}{c^2}A_2 - \ldots -\frac{p_n}{c^2}A_n.$$

Note that $\left(\mathcal{S}\left(\mathcal{A}\right),\odot_\psi\right)$ represents a commutative monoid. Moreover, the operator $\odot_\psi$ distributes over the sum $+_\psi$, that is, we have for all $A, B, C \in \mathcal{S}\left(\mathcal{A}\right)$, that $A \odot_\psi \left(B +_\psi C \right) = A \odot_\psi B +_\psi A \odot_\psi C$ \cite{laiate2021cross}. In fact,  $\left(\mathcal{S}\left(\mathcal{A}\right),+_\psi,\odot_\psi\right)$ yields a commutative ring.

Moreover, $\mathcal{S}\left(\mathcal{A}\right)$ is closed under the four arithmetic operations, that is, $B \otimes_\psi C \in \mathcal{S}\left(\mathcal{A}\right)$ for all $B, C \in \mathcal{S}\left(\mathcal{A}\right)$ and $\otimes_\psi \in \left\{+_\psi,-_\psi,\odot_\psi,\div_\psi\right\}$, called the set of the {\it $\psi-$arithmetic operations}. Recall that the $\psi-$arithmetic operations generalize the corresponding operations between real numbers
\cite{laiate2023properties}. Moreover, $\left(\mathcal{S}\left(\mathcal{A}\right),+_{\psi},\cdot_{\psi}\right)$ is an $n$-dimensional real vector space that has $\{A_1, \ldots, A_n\}$ as a basis \cite{esmi2021banach}.

Given any finite set of fuzzy numbers $\mathcal{B}=\left\{B_1,\ldots,B_n\right\} \subset \mathbb{R}_{\mathcal{F}}$ and $\varepsilon>0$ given, there exists an SLI set $\mathcal{A}=\left\{A_1,\ldots,A_n\right\}$ such that $\mathcal{D}\left(\mathcal{A},\mathcal{B}\right)< \varepsilon$, where $\mathcal{D}$ is a Hausdorff-based distance between finite sets of fuzzy numbers \cite{esmi2021banach}. In general, if no data set is given in a problem, we can generate an SLI set as follows:

\begin{theorem} \label{the:SLIgen}
    If $A \in \mathbb{R}_{\mathcal{F}}$ be a non-symmetric trapezoidal fuzzy number and $n \in \N$, then the sets
    \begin{equation} \label{eq:SLIgen}
        \left\{A^i\right\}_{i=1,\ldots,n} \quad \mbox{and} \quad \left\{\hat{f_i}\left(A\right)\right\}_{i=1,\ldots,n}
    \end{equation}
    are SLI, where 
    $A^{i}$ denotes the fuzzy modifier that defines the power hedges of $A$ for $i \in \mathbb{N}$.
\end{theorem}

The fuzzy number $A \in \mathbb{R}_{\mathcal{F}}$ of Theorem \ref{the:SLIgen} is also called {\it fuzzy basal number} when SLI sets are the underlying basis for a theory of calculus describing population dynamics (see \cite{laiate2023numerical} for details).
The next section introduces the notion of the linear fuzzy-valued function under the operation $\odot_{\psi}$.

\section{The Linear Fuzzy-valued Function under Cross Product}

Let $B,C \in \mathbb{R}_{\mathcal{F}}^{\wedge}$. The cross product between $B$ and $C$ is the fuzzy number given by
\begin{equation} \label{eq:crossprod}
    B \odot C = cB + bC -cb,
\end{equation}
\noindent where $\left[B\right]_1=\left\{b\right\}$ and $\left[C\right]_1=\left\{c\right\}$ \cite{ban2002cross,ban2006properties}. It is known that $\hat{f}\left(B,C\right) = B \odot C$, where $f(x,y)=cx+by-bc$, the linearization of the function $g(x,y)=xy$ around $(b,c)$ \cite{laiate2021bidimensional}. In this section, we extend this result to the general linear fuzzy-valued function $G: \mathbb{R}_{\mathcal{F}}^{\wedge} \rightarrow \mathbb{R}_{\mathcal{F}}^{\wedge}$ given by

\begin{equation} \label{eq:fuzzylinear}
    G(X) = A \odot X + C, \quad X \in \mathbb{R}_{\mathcal{F}}^{\wedge}.
\end{equation}

Let $P: \mathbb{R}^3 \to \mathbb{R}$ be the linearization of $g(x,y,z)=xy+z$ around $(x_0,y_0,z_0)$, that is, 
\begin{equation} \label{eq:linearization}
    \begin{split} 
        P(x,y,z) &= g(x_0,y_0,z_0) + \frac{\partial g}{\partial x}(x_0,y_0,z_0)(x-x_0)\\ 
        &+ \frac{\partial g}{\partial y}(x_0,y_0,z_0)(y-y_0) + \frac{\partial g}{\partial z}(x_0,y_0,z_0)(z-z_0)\\
        &= x_0y_0 + z_0 + y_0(x-x_0) + x_0(y-y_0) + (z-z_0)\\
        &= xy_0 +yx_0 - x_0y_0 + z
\end{split}
\end{equation}

The next proposition establishes a connection between the fuzzy-valued function given by Eq.\eqref{eq:fuzzylinear} and Eq. \eqref{eq:linearization}. 

\begin{proposition} \label{the:linearfunctionpsi}
    Let $X,A,B \in \mathbb{R}_{\mathcal{F}}^\wedge$ with $\left[X\right]_1=\left\{\overline{x}\right\}$, $\left[A\right]_1=\left\{\overline{a}\right\}$ and $\left[B\right]_1=\left\{\overline{b}\right\}$. The function $G(X) = A \odot X + B$ is the Zadeh extension of the linearization of $g(x,y,z)=xy +z$ around $(\overline{x},\overline{a},\overline{b})$, that is,
    \begin{equation} \label{eq:fuzzylinear2}
        G(X) = \hat{P}\left(X,A,B\right),
    \end{equation}
    where $P: \mathbb{R}^3 \to \mathbb{R}$ is given as in Equation \ref{eq:linearization} for $(x_0,y_0,z_0)=(\overline{x},\overline{a},\overline{b})$.
    \begin{proof}
        The proof follows similarly to Proposition 3 of \cite{laiate2021bidimensional}. In fact, the function $P(x,y,z)=\overline{a}x +\overline{x}y - \overline{a}\overline{x} + z$ is continuous in $\mathbb{R}^3$. By Theorem \ref{the:zadehsextension}, if $\hat{P}: \mathbb{R}_{\mathcal{F}} \times \mathbb{R}_{\mathcal{F}} \times \mathbb{R}_{\mathcal{F}} \to \mathbb{R}_{\mathcal{F}}$ is the Zadeh's extension of $P: \mathbb{R}^3 \to \mathbb{R}$, then
        \begin{equation}
            \hat{P}\left([X]_\alpha,[A]_\alpha,[B]_\alpha\right) = \overline{a}[X]_\alpha + \overline{x}[A]_\alpha - \overline{ax} + [B]_\alpha
        \end{equation}
        holds for all $\alpha \in [0,1]$, that is, $\hat{P}\left(X,A,B\right)=G(X)$.
    \end{proof}
\end{proposition}

Proposition \ref{the:linearfunctionpsi} suggests that the linear fuzzy-valued function under the cross product can be seen as a linearization of the linear fuzzy-valued function under the usual product $H: \mathbb{R}_{\mathcal{F}} \rightarrow \mathbb{R}_{\mathcal{F}}$, given by

\begin{equation} \label{eq:linearfunctionusual}
    H(X) = A\times X + B, \quad X \in \mathbb{R}_{\mathcal{F}}.
\end{equation}

When dealing with $\mathcal{S}\left(\mathcal{A}\right)-$linearly correlated fuzzy numbers, it is not possible to assume a direct product, since $\mathcal{S}\left(\mathcal{A}\right)$ is not closed under other types of arithmetic operations \cite{longo2022cross}. Due to the vector space structure of $\left(\mathcal{S}\left(\mathcal{A}\right),+_\psi,\cdot_\psi\right)$ a linearized product - the cross product - is necessary, and the operations involved shall be the $\psi-$arithmetic operations, which leads us to the $\psi-$cross product operation $\odot_\psi$.

Let $\mathcal{A} = \left\{1,A_2,\ldots,A_n\right\} \subset \mathbb{R}_{\mathcal{F}}^\wedge$ be SLI and $A, B \in \mathcal{S}\left(\mathcal{A}\right)$. Consider the function $F_{A,B}: \mathcal{S}\left(\mathcal{A}\right) \rightarrow \mathcal{S}\left(\mathcal{A}\right)$ given by

\begin{equation} \label{eq:fuzzylinearpsi}
    F_{A,B}(X) = A \odot_\psi X +_{\psi} B, \quad X \in \mathcal{S}\left(\mathcal{A}\right).
\end{equation}

In Eq. \eqref{eq:fuzzylinearpsi}, the fuzzy number $X \in \mathcal{S}\left(\mathcal{A}\right)$ can be seen as an independent variable of the function $F_{A,B}$. In addition, for each $X \in \mathcal{S}\left(\mathcal{A}\right)$, we can write 

\begin{equation} \label{eq:fuzzylinearpsi2}
    Y = A \odot_\psi X +_{\psi} B,
\end{equation}
\noindent from which the $\mathcal{S}\left(\mathcal{A}\right)-$linearly correlated fuzzy number $Y \in \mathcal{S}\left(\mathcal{A}\right)$ can be seen as a dependent variable of $X$. 

Since $A, B, X \in \mathcal{S}\left(\mathcal{A}\right)$, we have that $A = \psi\left(a_1,\ldots,a_n\right) = a_1 + a_2A_2 + \ldots a_nA_n$ and $B=\psi(b_1,\ldots,b_n) = b_1 + b_2A_2 + \ldots + b_nA_n$ for some $(a_1,\ldots,a_n), (b_1,\ldots,b_n) \in \mathbb{R}^n$. Moreover, for each $X \in \mathcal{S}\left(\mathcal{A}\right)$, there exists $(x_1,\ldots,x_n) \in \mathbb{R}^n$ such that $$X=\psi\left(x_1,\ldots,x_n\right) = x_1 + x_2A_2 + \ldots + x_nA_n.$$

Next, we prove the function $F_{A,B}$ is a bijection whenever $\left[A\right]_1 \neq \left\{0\right\}$.

\begin{proposition} \label{prop:bijection}
    Let $\mathcal{A} \subset \mathbb{R}_{\mathcal{F}}$ be SLI. The function $F_{A,B}: \Sc(\A) \to \Sc(\A) $ of Eq. \eqref{eq:fuzzylinearpsi}, is a bijection whenever $\left[A\right]_1 \neq \left\{0\right\}.$
    \begin{proof}
        First, let us prove the injectiveness of $F_{A,B}$. Let $X_1, X_2 \in \mathcal{S}\left(\mathcal{A}\right)$ be such that $F_{A,B}(X_1) = F_{A,B}(X_2)$. Note that 
        \begin{equation*} \label{eq:1}
            A \odot_\psi X_1 +_\psi B = A \odot_\psi X_2 +_\psi B,
        \end{equation*}
        Recall that $(\Sc(\A), +_\psi, \cdot_\psi)$ is a vector space by Proposition 1 which implies
        that $(\Sc(\A), +_\psi)$ is an abelian group. Therefore, 
        $A \odot_\psi X_1 = A \odot_\psi X_2$. Since $[A]_1 \neq \{0\}$, there exists $A_\psi^{-1} \in \Sc(\A)$ and we obtain
    \begin{eqnarray}
    \nonumber
        X_1&=&(A_{\psi}^{-1} \odot_\psi A) \odot_\psi X_1 = A_{\psi}^{-1} \odot_\psi (A \odot_\psi X_1) \\ \nonumber &=& A_{\psi}^{-1} \odot_\psi (A \odot_\psi X_1)=(A_{\psi}^{-1} \odot_\psi A) \odot_\psi X_2 = X_2.
         \end{eqnarray}
        Therefore, $F_{A,B}$ is one-to-one. To prove that $F_{A,B}$ is onto, consider an arbitrary element  $Y$ of $\mathcal{S}\left(\mathcal{A}\right)$. It suffices to show that $Y = F_{A,B}(X)$ for some $X \in \mathcal{S}\left(\mathcal{A}\right)$. Consider $X = A_{\psi}^{-1}\odot_\psi Y -_\psi A_{\psi}^{-1}\odot_\psi B$. It follows from observations in Section 2 that $X \in \mathcal{S}\left(\mathcal{A}\right)$. In addition, we have
        \begin{equation*}
        \begin{split}
            F_{A,B}(X) &= F_{A,B}(A_{\psi}^{-1}\odot_\psi Y -_{\psi} A_{\psi}^{-1}\odot_\psi B)\\ 
            &= A \odot_\psi \left(A_{\psi}^{-1}\odot_\psi Y - A_{\psi}^{-1}\odot_\psi B\right) +_\psi B\\ 
            &= (A_\psi^{-1} \odot_\psi A) \odot_\psi Y -_\psi (A_\psi^{-1} \odot_\psi A) \odot_\psi B +_\psi B = Y,
        \end{split} 
        \end{equation*}
\noindent from which we conclude $F_{A,B}$ is surjective. Therefore, $F_{A,B}$ is a bijection on $\mathcal{S}\left(\mathcal{A}\right)$.
    \end{proof}
\end{proposition}

Proposition \ref{prop:bijection} shows that the linear fuzzy-valued function $F_{A,B}$ establishes a bijection in $\mathcal{S}\left(\mathcal{A}\right)$ whenever $\mathcal{A} \subset \mathbb{R}_{\mathcal{F}}$ is SLI and $\left[A\right]_1 \neq \left\{0\right\}$. Moreover, the inverse function $F_{A,B}^{-1}: \mathcal{S}\left(\mathcal{A}\right) \rightarrow \mathcal{S}\left(\mathcal{A}\right)$ is given by
\begin{equation} \label{eq:inverseF}
    F_{A,B}^{-1}(X) = A_\psi^{-1} \odot_\psi X -_\psi A_{\psi}^{-1} \odot_\psi B, \quad \forall X \in \mathcal{S}\left(\mathcal{A}\right).
\end{equation}
The next examples evaluate $F_{A,B}\left(X_\lambda\right)$, where $\left\{X_\lambda\right\}$ is given by two different indexed families of $\Sc(\A)-$linearly correlated fuzzy numbers.

\begin{example} \label{ex:example1}
    Let $X \in \mathcal{S}\left(\mathcal{A}\right)$ be given and consider the family of $\mathcal{S}-$linearly correlated fuzzy numbers given by $\left\{X_\lambda\right\}_{\lambda \in \mathbb{R}} = \left\{\lambda X : \lambda \in \mathbb{R}\right\}$. Then,
    $$
    F_{A,B}\left(X_\lambda\right) = A \odot_\psi \left(X_\lambda\right) +_\psi B = A \odot_\psi \left(\lambda X\right) +_\psi B = \lambda \left(A \odot_\psi X\right) +_\psi B, \quad \forall \lambda \in \mathbb{R}.
    $$
    In this case, we identify the function $F_{A,B}$ with the fuzzy number-valued application $f: \mathbb{R} \to \mathcal{S}\left(\mathcal{A}\right)$ given by $f(t)=\left(A \odot_\psi X\right)t +_\psi B$, $\forall t \in \mathbb{R}$.
    


\end{example} 

\begin{example} \label{ex:example2}
    Let $X \in \mathcal{S}\left(\mathcal{A}\right)$ be given and consider the family of $\mathcal{S}-$linearly correlated fuzzy numbers given by $\left\{X_\mu\right\}_{\mu \in \mathbb{R}} = \left\{X +_{\psi} \mu : \mu \in \mathbb{R}\right\}$. Note that since $\mathcal{A}=\left\{1,A_2,\ldots,A_n\right\}$, $\mu = \psi\left(\mu,0,\ldots,0\right) \in \mathcal{S}\left(\mathcal{A}\right)$. In this case, we have that
    $$
    F_{A,B}\left(X_\mu\right) = A \odot_\psi \left(X_\mu\right) +_\psi B = A \odot_\psi \left(X +_\psi \mu\right) +_\psi B = \left(A \odot_\psi X +_\psi B\right) +_\psi \mu A, \quad \forall \mu \in \mathbb{R}.
    $$
    In this case, we identify the function $F_{A,B}$ with the fuzzy number-valued application $f: \mathbb{R} \to \mathcal{S}\left(\mathcal{A}\right)$ given by $f(t)=\left(A \odot_\psi X+_\psi B\right) +_\psi At$, $\forall t \in \mathbb{R}$. 
    

\end{example}

\section{Fuzzy Arithmetic Equations in $\mathcal{S}(\mathcal{A})$}

Let $\mathcal{A} \subset \mathbb{R}_{\mathcal{F}}$ be a given SLI set. Consider the fuzzy arithmetic equation

\begin{equation} \label{eq:fuzarithmeq}
    A \odot_\psi X +_\psi B = C,
\end{equation}
\noindent where $A, B, C \in \mathcal{S}\left(\mathcal{A}\right)$ are given. If $\left[A\right]_1 \neq \left\{0\right\}$, Eq. \eqref{eq:fuzarithmeq} has only one solution, given by
\begin{equation} \label{eq:fuzarithmeq2}
    X=A_\psi^{-1} \odot_\psi \left(C -_\psi B\right) \, \cite{laiate2021cross}.    
\end{equation}
Denote $\overline{B} = C -_\psi B$. Note that solving Eq. \eqref{eq:fuzarithmeq} is equivalent to solving the fuzzy equation 
\begin{equation} \label{eq:fuzarithmeq3}
    A \odot_\psi X = \overline{B}.
\end{equation}
The solution to the Eq. \eqref{eq:fuzarithmeq} or, equivalently, the solution to the Eq. \eqref{eq:fuzarithmeq3} is the inverse image of $C \in \mathcal{S}\left(\mathcal{A}\right)$ via $F_{A,B}: \mathcal{S}\left(\mathcal{A}\right) \to \mathcal{S}\left(\mathcal{A}\right)$, given by Eq. \eqref{eq:fuzzylinearpsi}. Thus, the solution, when exists, is given by $X = F_{A,B}^{-1}\left(C\right)$. By Proposition \ref{prop:bijection}, the existence and uniqueness of the solution is assured whenever the condition $\left[A\right]_1 \neq \left\{0\right\}$ is satisfied.

On the other hand, if $\left[A\right]_1 = \left\{0\right\}$, then $$A \odot_\psi X = \left[A\right]_1X +_\psi \left[X\right]_1A - \left[A\right]_1\left[X\right]_1 = Ax,$$
so that Eq. \eqref{eq:fuzzylinear2} can be written as
\begin{equation} \label{eq:fuzarithmeq3}
    Ax=\overline{B},
\end{equation}
\noindent which \eqref{eq:fuzarithmeq3} has a solution if and only if there exists $x \in \mathbb{R}$ satisfying Eq. \eqref{eq:fuzarithmeq3}. This implies that $A$ and $\overline{B}$ are linearly correlated, that is, $\overline{B} \in \mathbb{R}_{\mathcal{F}(A)}$, where
\begin{equation} \label{Rfa}
    \mathbb{R}_{\mathcal{F}(A)} = \left\{B = qA + r : q, r \in \mathbb{R}\right\} = \mathcal{S}\left(A,1\right),
\end{equation}
\noindent the set of all linear combinations of $A$ and $1$, where $1 \in \mathbb{R}_{\mathcal{F}}$ is regarded as a singleton \cite{esmi2018frechet}.

\subsection*{Further Discussion}

Denoting $[X]_1 = \left\{\overline{x}\right\}$, $[A]_1 = \left\{\overline{a}\right\}$, and assuming that $A_1 =1$, Eq. \eqref{eq:fuzarithmeq3} can be written levelwise as
\begin{equation} \label{eq:fuzaritheq4}
    \sum_{i=1}^{n} \left(\overline{x}a_i + \overline{a}x_i\right) \left[A_i\right]_{\alpha} - \overline{ax} = [B_i]_{\alpha}, \quad \forall \alpha \in [0,1].
\end{equation}
Eq. \eqref{eq:fuzaritheq4} and Proposition \ref{the:linearfunctionpsi} show that solving a linear fuzzy arithmetic equation in $\mathcal{S}\left(\mathcal{A}\right)$ is equivalent to solving a family of linear interval equation with real coefficients. Although Eq. \eqref{eq:fuzarithmeq3} involves an operation of product between fuzzy numbers, it does not involve a product between intervals. This reasoning can be extended to fully fuzzy linear systems in $\mathcal{S}\left(\mathcal{A}\right)$ under the $\psi-$operations, since Eq. \eqref{eq:fuzarithmeq3} can be seen as a fuzzy linear system of dimension $1$. 

\section{Concluding Remarks}

This paper investigates linear equations given by Eq. \eqref{eq:main} in the commutative ring $(\Sc(\A), +_\psi, \odot_\psi)$, where
$\A$ is an SLI set of the form  $\left\{1,A_2,\ldots,A_n\right\} \subset \R_{\F}^\wedge$
\cite{esmi2021banach,esmi2022calculus,laiate2021cross}. 
To this end, we defined a function $F_{A,B}: \Sc(\A) \to \Sc(\A)$ for any $A, B \in \Sc(\A) \subset  \R_{\F}^\wedge$ given by 
$F_{A,B}(X)=A \odot_\psi X +_\psi B$ for all $X \in \mathcal{S}\left(\mathcal{A}\right)$
and proved that $F_{A,B}$ is bijective whenever the core of $A \in \mathcal{S}\left(\mathcal{A}\right)$ 
differs from $\{ 0 \}$. Additionally, we extended the result presented in \cite{laiate2021bidimensional} by proving that, when considering the cross product $\odot$ instead of the $\psi-$cross product $\odot_\psi$, the function $F_{A,B}$ equals Zadeh's extension of the function $f(x,y,z)=xy+z$ around $(\overline{x},\overline{a},\overline{b})$, where $[X]_1=\{\overline{x}\}, [A]_1=\{\overline{a}\}$, and $[B]_1=\{\overline{b}\}$.

Lastly, we described the solution to the linear fuzzy arithmetic equation given by Eq. \eqref{eq:main}
for $A \in \R_{\F}^\wedge$ with $[A]_1 \neq \{0\}$
in the $(\Sc(\A), +_\psi, \odot_\psi)$ setting  as  $F^{-1}_{A,B}(C)$ for which we provided a formula. It is worth noting that Eq. \eqref{eq:main} does not consist of a direct generalization of an interval linear equation with interval coefficients \cite{lodwick2015interval}. Since the cross product is a linearized operation, the $\alpha-$levels of Eq. \eqref{eq:main} are given by linear combinations of intervals with real coefficients, that is, they do not include products of intervals.

In future research, we intend to extend the results of this paper. In particular, we plan to address systems of fuzzy linear equations having values in $\Sc(\A)$.

\begin{acknowledgement}
The authors would like to thank the School of Applied Mathematics at the Getúlio Vargas Foundation (FGV/EMAp), the Institute of Mathematics, Statistics and Scientific Computing of the University of Campinas (IMECC/Unicamp), as well as Prof. Dr. Fernando Gomide for inspiring the contents of this paper.
This work was supported in part by FAPESP under grant no. 2020/09838–0 (Brazilian Institute of Data Science).
\end{acknowledgement}

\bibliography{splncs04}    

\begin{thebibliography}{10}

\bibitem{ban2002cross}
{\sc Ban, A., and Bede, B.}
\newblock Cross product of l-r fuzzy numbers and applications.
\newblock {\em Anal. Univ. Oradea, fasc. math 9\/} (2002), 95--108.

\bibitem{ban2006properties}
{\sc Ban, A., and Bede, B.}
\newblock Properties of the cross product of fuzzy numbers.
\newblock {\em Journal of Fuzzy Mathematics 14}, 3 (2006).

\bibitem{bede}
{\sc Bede, B.}
\newblock {\em Mathematics of Fuzzy Sets and Fuzzy Logic}.
\newblock Studies in Fuzziness and Soft Computing. Springer,Berlin, Heidelberg,
  2013.

\bibitem{behera2014solving}
{\sc Behera, D., and Chakraverty, S.}
\newblock Solving fuzzy complex system of linear equations.
\newblock {\em Information Sciences 277\/} (2014), 154--162.

\bibitem{biacino1989equations}
{\sc Biacino, L., and Lettieri, A.}
\newblock Equations with fuzzy numbers.
\newblock {\em Information Sciences 47}, 1 (1989), 63--76.

\bibitem{buckley1990solving}
{\sc Buckley, J., and Qu, Y.}
\newblock Solving linear and quadratic fuzzy equations.
\newblock {\em Fuzzy Sets and Systems 38}, 1 (1990), 43--59.

\bibitem{diamond1990metric}
{\sc Diamond, P., and Kloeden, P.}
\newblock Metric spaces of fuzzy sets.
\newblock {\em Fuzzy Sets and Systems 35}, 2 (1990), 241--249.

\bibitem{dubois1978operations}
{\sc Dubois, D., and Prade, H.}
\newblock Operations on fuzzy numbers.
\newblock {\em International Journal of Systems Science 9}, 6 (1978), 613--626.

\bibitem{dubois1983inverse}
{\sc Dubois, D., and Prade, H.}
\newblock Inverse operations for fuzzy numbers.
\newblock {\em IFAC Proceedings Volumes 16}, 13 (1983), 399--404.

\bibitem{dubois1993fuzzy}
{\sc Dubois, D., and Prade, H.}
\newblock Fuzzy numbers: an overview.
\newblock {\em Readings in Fuzzy Sets for Intelligent Systems\/} (1993),
  112--148.

\bibitem{esmi2021banach}
{\sc Esmi, E., de~Barros, L.~C., Santo~Pedro, F., and Laiate, B.}
\newblock Banach spaces generated by strongly linearly independent fuzzy
  numbers.
\newblock {\em Fuzzy Sets and Systems 417\/} (2021), 110--129.

\bibitem{esmi2019some}
{\sc Esmi, E., de~Barros, L.~C., and Wasques, V.~F.}
\newblock Some notes on the addition of interactive fuzzy numbers.
\newblock In {\em International Fuzzy Systems Association World Congress\/}
  (2019), Springer, pp.~246--257.

\bibitem{esmi2022calculus}
{\sc Esmi, E., Laiate, B., Santo~Pedro, F., and Barros, L.~C.}
\newblock Calculus for fuzzy functions with strongly linearly independent fuzzy
  coefficients.
\newblock {\em Fuzzy Sets and Systems 436\/} (2022), 1--31.

\bibitem{esmi2018frechet}
{\sc Esmi, E., Santo~Pedro, F., de~Barros, L.~C., and Lodwick, W.}
\newblock Fr{\'e}chet derivative for linearly correlated fuzzy function.
\newblock {\em Information Sciences 435\/} (2018), 150--160.

\bibitem{fuller1998fuzzy}
{\sc Full{\'e}r, R., et~al.}
\newblock {\em Fuzzy reasoning and fuzzy optimization}.
\newblock No.~9. Turku Centre for Computer Science Abo, 1998.

\bibitem{ghanbari2022new}
{\sc Ghanbari, M., Allahviranloo, T., Nuraei, R., and Pedrycz, W.}
\newblock A new effective approximate multiplication operation on lr fuzzy
  numbers and its application.
\newblock {\em Soft Computing 26}, 9 (2022), 4103--4113.

\bibitem{laiate2023properties}
{\sc Laiate, B.}
\newblock On the properties of fuzzy differential equations under cross
  operations.
\newblock {\em Computational and Applied Mathematics 42}, 6 (2023).

\bibitem{laiate2021bidimensional}
{\sc Laiate, B., Barros, L.~C., Santo~Pedro, F., and Esmi, E.}
\newblock Bidimensional fuzzy initial value problem of autocorrelated fuzzy
  processes via cross product: the prey-predator model.
\newblock In {\em Proceedings of the 19th World Congress of the International
  Fuzzy Systems Association (IFSA), 12th Conference of the European Society for
  Fuzzy Logic and Technology (EUSFLAT), and 11th International Summer School on
  Aggregation Operators (AGOP)\/} (2021), Atlantis Press, pp.~171--178.

\bibitem{laiate2023numerical}
{\sc Laiate, B., Longo, F., Alves, J.~R., and Meyer, J. F.~C.}
\newblock Numerical solutions of fuzzy population models: A case study for
  chagas’ disease dynamics.
\newblock In {\em Fuzzy Information Processing 2023: Proceedings of the North
  American Fuzzy Information Processing Society Annual Conference\/} (2023),
  Springer, pp.~172--183.

\bibitem{laiate2021cross}
{\sc Laiate, B., Watanabe, R.~A., Esmi, E., Santo~Pedro, F., and Barros, L.~C.}
\newblock A cross product of s-linearly correlated fuzzy numbers.
\newblock In {\em 2021 IEEE International Conference on Fuzzy Systems
  (FUZZ-IEEE)\/} (2021), IEEE, pp.~1--6.

\bibitem{lodwick2015interval}
{\sc Lodwick, W.~A., and Dubois, D.}
\newblock Interval linear systems as a necessary step in fuzzy linear systems.
\newblock {\em Fuzzy Sets and Systems 281\/} (2015), 227--251.

\bibitem{longo2022cross}
{\sc Longo, F., Laiate, B., Pedro, F.~S., Esmi, E., Barros, L.~C., and Meyer,
  J.~F.}
\newblock A-cross product for autocorrelated fuzzy processes: the hutchinson
  equation.
\newblock In {\em Explainable AI and Other Applications of Fuzzy Techniques:
  Proceedings of the 2021 Annual Conference of the North American Fuzzy
  Information Processing Society, NAFIPS 2021\/} (2022), Springer,
  pp.~241--252.

\bibitem{nguyen1978note}
{\sc Nguyen, H.~T.}
\newblock A note on the extension principle for fuzzy sets.
\newblock {\em Journal of Mathematical Analysis and Applications 64}, 2 (1978),
  369--380.

\bibitem{sanchez1984solution}
{\sc Sanchez, E.}
\newblock Solution of fuzzy equations with extended operations.
\newblock {\em Fuzzy Sets and Systems 12}, 3 (1984), 237--248.

\bibitem{sevastjanov2009new}
{\sc Sevastjanov, P., and Dymova, L.}
\newblock A new method for solving interval and fuzzy equations: linear case.
\newblock {\em Information Sciences 179}, 7 (2009), 925--937.

\bibitem{wang09}
{\sc Wang, X., Ruan, D., and Kerre, E.}
\newblock {\em Mathematics of Fuzziness—Basic Issues}.
\newblock Studies in Fuzziness and Soft Computing. Springer Berlin Heidelberg,
  2009.

\bibitem{yager1980lack}
{\sc Yager, R.~R.}
\newblock On the lack of inverses in fuzzy arithmetic.
\newblock {\em Fuzzy Sets and Systems 4}, 1 (1980), 73--82.

\bibitem{zadeh1975concept}
{\sc Zadeh, L.~A.}
\newblock The concept of a linguistic variable and its application to
  approximate reasoning—i.
\newblock {\em Information Sciences 8}, 3 (1975), 199--249.

\end{thebibliography}

\end{document}